 \documentclass[10pt]{amsart}
\usepackage{mathrsfs}
\usepackage{amsfonts}
\usepackage{amsmath}
\usepackage{amssymb}
\usepackage{amsthm}
\usepackage{enumerate}
\usepackage{color}
\usepackage{geometry}
\usepackage{hyperref}
\usepackage{float}
\usepackage{graphicx}
\usepackage{multirow}
\usepackage[numbers,sort&compress]{natbib}
\usepackage{color}
\usepackage[all]{xy}
\usepackage{cases}
\usepackage{multirow}

\allowdisplaybreaks
\numberwithin{equation}{section}

\theoremstyle{plain}
\newtheorem{prop}{Proposition}[section]
\newtheorem{coro}[prop]{Corollary}

\newtheorem{lemm}[prop]{Lemma}

\newtheorem{theorem}[prop]{Theorem}

\theoremstyle{definition}
\newtheorem{definition}[prop]{Definition}

\title{A construction of a free digroup$^*$}

\author{Guangliang Zhang}
\address{G.Z., School of Mathematics and Systems Science, Guangdong Polytechnic Normal University, Guangzhou 510631, P. R. China}
\email{zgl541@163.com}

\author{ Yuqun Chen$^{\sharp}$}
\address{Y.C., School of Mathematical Sciences, South China Normal University, Guangzhou 510631, P. R. China}
\email{yqchen@scnu.edu.cn}

\thanks{${}^*$ Supported by the NNSF of China (11571121), the NSF of Guangdong Province (2017A030313002) and the Science and Technology Program of Guangzhou (201707010137)}

\thanks{${}^{\sharp}$ Corresponding author}



\begin{document}

\maketitle

$$
\mbox{Communicated by  Mikhail Volkov}
$$

\begin{abstract}
We give a construction of a free digroup on a set $X$ and formulate the halo and the group parts of it. We prove that a free digroup on $X$ is isomorphic to a free digroup on a set $Y$ if and only if $card(X)=card(Y)$.
\end{abstract}

%%%%

\section{Introduction and Preliminaries}

In the theory of Leibniz algebras,
one of the prominent open problems is to find an appropriate generalization of Lie’s third theorem,
which associates a (local) Lie group to any (real or complex) Lie algebra.
A key and difficult aspect of this problem is to find the appropriate analogue of Lie group for Leibniz algebras,
that is, to determine what should be the correct generalization of the notion of a group.
So little is known about what properties these group-like objects should have that Loday dubbed them ``coquecigrues" in \cite{Lo93}.
And so the problem has come to be known as the ``coquecigrue" problem (for Leibniz algebras).

The notion of a digroup first implicitly appeared in Loday's work \cite{Lo99}.
Further, it was proposed independently by Kinyon \cite{Kinyon04}, Felipe \cite{Felipe} and Liu \cite{Liu},
to provide a partial solution to the coquecigrue problem for Leibniz algebras.
Kinyon \cite{Kinyon04} gave a much clearer definition of a digroup as follow.
\begin{definition}
A \textit{digroup} is a pair $(G,1)$ equipped with two binary operations $\vdash$ and $\dashv$, a unary operation $\dag$, and a nullary operation 1, where~$1$ is called the {\it unit} of $G$, satisfying each of the following six axioms:
\begin{equation}\label{eq00}
\begin{cases}
G1. \ (G, \vdash)\  \mbox{and} \  (G, \dashv) \  \mbox{are both semigroups}, \\
G2. \ a\dashv(b\vdash c)=a\dashv (b\dashv c),\\
G3. \ (a\dashv b)\vdash c=(a\vdash b)\vdash c, \\
G4. \ a\vdash(b\dashv c)=(a\vdash b)\dashv c, \\
G5. \ 1 \vdash a = a = a \dashv 1, \\
G6. \ a \vdash a^{\dag} = 1 = a^{\dag} \dashv a.
\end{cases}
\end{equation}
\end{definition}
Using semigroup theory, Kinyon showed that every digroup is a product of a group and a ``trivial" digroup,
in the service of his partial solution to the coquecigrue problem.
Of course, digroups are a generalization of groups and play an important role in this open problem from the theory of Leibniz algebras.
Recently the study of algebraic properties on digroups  has attracted considerable attention.
In \cite{Ph}, Phillips gave a simple basis of independent axioms for the variety of digroups.
Salazar-D\'\i az, Vel\'asquez and Wills-Toro \cite{Salazar} studied a further generalization of the digroup structure
and showed analogues to the first isomorphism theorem.
Lately Ongay, Vel\'asquez and Wills-Toro \cite{Ongay} discussed the notion of normal subdigroups and quotient digroups,
and then  established the corresponding analogues of the classical Isomorphism Theorems.
Different examples of digroups can be found in \cite{AYZhuchok17}.
For further investigation of the structure of digroups,
A.V. Zhuchok proposed some open problems at the end of his paper \cite{Zhu17}.
One of the open problems is: ``Construct a free digroup".

In this paper, we apply the method of Gr\"{o}bner-Shirshov bases for dialgebras in \cite{CZ} to
give a construction of a free digroup on a set
and solve the above Zhuchok's problem.

In \cite{CZ}, a Composition-Diamond lemma for dialgebras was established and
a method to find normal forms of elements of an arbitrary disemigroup was given.
Clearly, the variety of digroups is a subvariety of the variety of disemigroups.
Then the free digroup on a set is a quotient of a free disemigroup  modulo certain relations.

The paper is organized as follows. In section 2, we recall the definitions of the Gr\"{o}bner-Shirshov bases
and Composition-Diamond Lemma for dialgebras in \cite{CZ}.
In section 3, we give a Gr\"{o}bner-Shirshov basis for the free digroup $F(X)$ on a set $X$ and thus obtain  explicit normal forms of elements of $F(X)$.
Moreover, we also formulate the halo and the group parts of $F(X)$
and show that $F(X)\cong F(Y)$ if and only if $card(X)=card(Y)$.

\section{Composition-Diamond lemma for dialgebras}

As is known to all, the method of Gr\"{o}bner-Shirshov bases is a powerful
tool to solve the normal form problem in various categories, including dialgebras \cite{CZ}.
In this section, we review some of conceptions and results on Gr\"{o}bner-Shirshov bases for dialgebras, see  \cite{CZ}.

Let $\mathbf{k}$ be a field. A \textit{dialgebra} (\textit{disemigroup}) is a $\mathbf{k}$-module (set) equipped with two binary operations $\vdash$ and $\dashv$, satisfying axioms $G1$-$G4$ in (\ref{eq00}). Recall that for every dialgebra (disemigroup)~$D$, for all~$x_1,...,x_n$ in~$D$, every parenthesizing of
$
x_1\vdash\cdots\vdash x_m\dashv\cdots\dashv x_t
$
 gives the same element in~$D$, which we denote by~$[x_1...x_t]_m$. In particular, the notation~$[x_1]_1$ means~$x_1$.
Let $Di\langle X\rangle$ be the free dialgebra over $\mathbf{k}$ generated by a set $X$,
$X^+$ the set of all nonempty associative words on $X$ and $X^*=X^+\cup \{\epsilon\}$, the free monoid on $X$, where $\epsilon$ is the empty word. Write
$$
[X^+]_\omega:=\{[u]_m \mid u\in X^+, m\in \mathbb{Z}^+, 1\leq m \leq |u|\},
$$
where $|u|$ is the length of $u$. Note that $[x]_1=x$ if $x\in X$. It is well known from \cite{Lo99,CZ} that $[X^+]_\omega$ is the free disemigroup on $X$ and a $\mathbf{k}$-basis of $Di\langle X\rangle$, where for all $[u]_m, [v]_n \in [X^+]_\omega$,
$$
[u]_m\vdash [v]_n=[uv]_{|u|+n}, \ \ \ [u]_m\dashv [v]_n=[uv]_m.
$$
Free disemigroups were studied in  \cite{Zhu11}.
For any $h=[u]_m \in [X^+]_\omega$,
we call $u$ the \textit{associative word} of $h$, and
$m$, denoted by $p(h)$, the \textit{position of center} of $h$.

Let $X$ be a well-ordered set.
For any $u=x_{i_1}x_{i_2}\cdots x_{i_n}, v=x_{j_1}x_{j_2}\cdots x_{j_m}\in X^+$, where  $x_{i_l},x_{j_t}\in X$, define
\begin{equation*}\label{equ0}
u>v \ \Leftrightarrow \ (|u|,x_{i_1},x_{i_2},\cdots, x_{i_n})>(|v|,x_{j_1},x_{j_2},\cdots, x_{j_m}) \ \mbox{lexicographically}.
\end{equation*}
Then $>$ is a well ordering on $X^+$ and we call it the \textit{deg-lex ordering}.

A well ordering $>$ on $X^+$ is \textit{monomial} if for any $u,v\in X^+$, we have
$$
u>v \Rightarrow w_1uw_2>w_1vw_2 \ \mbox{for all} \ w_1,w_2\in X^*.
$$
Clearly, the deg-lex ordering is monomial.
Let $X^+$ be endowed with a monomial ordering. We define the \textit{monomial-center ordering} on $[X^+]_\omega$ as follow: for any $[u]_m,[v]_n\in [X^+]_\omega$,
\begin{equation}\label{equ0}
[u]_m>[v]_n \ \Leftrightarrow \  (u,m)>(v,n) \ \ \mbox{lexicographically}.
\end{equation}
In particular, if $X^+$ is endowed with the deg-lex ordering,
we call the ordering defined by $(\ref{equ0})$ the \textit{deg-lex-center ordering} on $[X^+]_\omega$.
It is clear that a monomial-center ordering is a well ordering on $[X^+]_\omega$.
Here and subsequently, the monomial-center ordering on $[X^+]_\omega$ will be used, unless otherwise stated.

For convenience we assume that $[u]_m>0$ for any $[u]_m\in [X^+]_\omega$.
Then every nonzero polynomial $f\in Di\langle X\rangle$ has the leading monomial $\overline{f}$.
We denote the associative word of $\overline{f}$ by $\widetilde{f}$ and the leading term of $f$ by $lt(f)$.
If $\widetilde{f}> \widetilde{r\!_{_f}}$, where $r\!_{_f}:=f-lt(f)$, then $f$ is called \textit{strong}.
If the coefficient of $\overline{f}$ in $f$ is equal to 1, then $f$ is called \textit{monic}.
For a subset $S$ of $Di\langle X\rangle$,
$S$ is\textit{ monic} if $s$ is monic for all $s\in S$.

Let $s$ be a monic polynomial. If
\begin{eqnarray}\label{d1}
g=[x_{i_1}\cdots x_{i_k}\cdots x_{i_n}]_m|_{_{x_{i_k}\mapsto s}}:=[x_{i_1}\cdots x_{i_{k-1}}sx_{i_{k+1}}\cdots x_{i_n}]_m,
\end{eqnarray}
where $1\leq k \leq n,\ x_{i_l}\in X,\ 1\leq l \leq n$, then we call $g$ an \textit{$s$-diword}. For simplicity, we denote the $s$-diword of the form (\ref{d1}) by $(asb)$, where $a,b\in X^*$. Moreover, if either $k=m$ or $s$ is strong, then we call the $s$-diword $g$ is \textit{normal}.
Let $(asb)$ be a normal $s$-diword, then $\overline{(asb)}=[a\widetilde{s}b]_l$ for some $l\in P([asb])$, where
\begin{displaymath} P([asb]):=
\begin{cases}
\{n\in \mathbb{Z}^+\mid 1\leq n \leq |a|\} \cup \{|a|+p(\overline{s})\} \cup \{n\in \mathbb{Z}^+\mid |a\widetilde{s}|<n\leq |a\widetilde{s}b|\}
 \ \ \text{if $s$ is strong,}\\
\{|a|+p(\overline{s})\}  \ \  \text{if $s$ is not strong.}
\end{cases}
\end{displaymath}
If this is so, we denote the normal $s$-diword (or normal $S$-diword) $(asb)$ by $[asb]_l$.

Let $f,g$ be monic polynomials in $Di\langle X\rangle$.
\begin{enumerate}
\item[1)] If $f$ is not strong, then we call $x \dashv f$ the \textit{composition of left multiplication} of $f$
for all $x\in X$ and $f \vdash [u]_{|u|}$ the \textit{composition of right multiplication} of $f$  for all $u\in X^+$.

\item[2)] Suppose that  $w=\widetilde{f}= a\widetilde{g}b$ for some $a,b\in X^*$
and $(agb)$ is a normal $g$-diword.
\begin{enumerate}
\item[2.1] If $p(\overline{f})\in P([agb])$,
then we call
$$
(f,g)_{\overline{f}}=f-[agb]_{p(\overline{f})}
$$
the \textit{composition of inclusion} of $f$ and $g$.

\item[2.2] If $p(\overline{f})\notin P([agb])$ and both $f$ and $g$ are strong,
then for any $x\in X$ we call
$$
(f,g)_{[xw]_1}=[xf]_1-[xagb]_1
$$
the \textit{composition of left multiplicative inclusion} of $f$ and $g$, and
$$
(f,g)_{[wx]_{_{|wx|}}}=[fx]_{|wx|}-[agbx]_{|wx|}
$$
the \textit{composition of right multiplicative inclusion} of $f$ and $g$.
\end{enumerate}

\item[3)] Suppose that there exists $w=\widetilde{f}b= a\widetilde{g}$ for some $a,b\in X^*$
such that $|\widetilde{f}|+|\widetilde{g}|>|w|$,  $(fb)$ is a normal $f$-diword
and $(ag)$ is a normal $g$-diword.
\begin{enumerate}
\item[3.1]
 If $P([fb])\cap P([ag])\neq \varnothing$,
then for any $m\in P([fb])\cap P([ag])$ we call
$$
(f,g)_{[w]_m}=[fb]_m-[ag]_m
$$
the \textit{composition of intersection} of $f$ and $g$.

\item[3.2] If $P([fb])\cap P([ag])=\varnothing$ and both $f$ and $g$ are strong,
then for any $x\in X$ we call
$$
(f,g)_{[xw]_1}=[xfb]_1-[xag]_1
$$
the \textit{composition of left multiplicative intersection} of $f$ and $g$, and
$$
(f,g)_{{[wx]}_{|wx|}}=[fbx]_{|wx|}-[agx]_{|wx|}
$$
the \textit{composition of right multiplicative intersection} of $f$ and $g$.
\end{enumerate}
\end{enumerate}

For any composition $(f,g)_{[u]_n}$ mentioned above, we call $[u]_n$ the \textit{ambiguity} of $f$ and $g$.

\begin{definition}\label{dcgsb}\cite{CZ}
Let $S$ be  a monic subset of $Di\langle X\rangle$.
A polynomial $h\in Di\langle X\rangle$ is called \emph{trivial modulo} $S$,
if $h=\sum_i \alpha_{i}[a_i s_i b_i]_{m_i}$, where each $\alpha_i\in \mathbf{k}, \ a_i,b_i\in X^*, \ s_i\in S$,
and $\overline{[a_i s_i b_i]_{m_i}}\leq \overline{h}$ if $\alpha_{i}\neq0$.

A monic set $S$ is called a \emph{Gr\"{o}bner-Shirshov basis} in $Di\langle X\rangle$ if any
composition of polynomials in $S$ is trivial modulo $S$.

\end{definition}

For convenience, for any $f,g\in Di\langle X\rangle$ and $[w]_m\in [X^+]_\omega$, we write
$$
f\equiv g  \mod(S,[w]_m)
$$
which means that $f-g=\sum_i \alpha_{i}[a_i s_i b_i]_{m_i}$, where each $\alpha_i\in \mathbf{k}, \ a_i,b_i\in X^*, \ s_i\in S$,
and $\overline{[a_i s_i b_i]_{m_i}}< [w]_m\ $ if $\alpha_{i}\neq0$.

\ \

Note that for a monic set $S$, $S$ is a Gr\"{o}bner-Shirshov basis in $Di\langle X\rangle$ if and only if
\begin{enumerate}
\item[(i)] \ for any composition $h$ of the form 1) in $S$, $h$ is trivial modulo $S$;
\item[(ii)] \ for any $f,g\in S$ and any composition $(f,g)_{[w]_m}$ of the form 2) or 3), $(f,g)_{[w]_m}\equiv 0  \mod(S,[w]_m)$.
\end{enumerate}

\ \

The following lemma is Composition-Diamond lemma for dialgebras in \cite{CZ}.

\begin{lemm}(\cite[Theorem 3.18]{CZ})\label{cd}
Let $S$ be a monic subset of $Di\langle X\rangle$,
$>$ a monomial-center ordering on $[X^+]_\omega$ and
$Id(S)$ the ideal of $Di\langle X\rangle$ generated by $S$. Then the following statements are equivalent.
\begin{enumerate}
\item[(i)] \ $S$ is a Gr\"{o}bner-Shirshov basis in $Di\langle X\rangle$.
\item[(ii)] \ $f\in Id(S)\Rightarrow \overline{f}=\overline{[asb]_m}$ for some normal $S$-diword $[asb]_m$.
\item[(iii)] \
$Irr(S)=\{[u]_n\in [X^+]_\omega\mid [u]_n\neq \overline{[asb]_m} \ \mbox{ for any normal } S\mbox{-diword}\ [asb]_m \}$
is a $\mathbf{k}$-basis of the quotient dialgebra $Di\langle X\mid  S\rangle:=Di\langle X\rangle/Id(S)$.
\end{enumerate}
\end{lemm}

Let $Disgp\langle X\rangle:=[X^+]_\omega$ be the free disemigroup on $X$.
For an arbitrary disemigroup $D$, $D$ has an expression
$$
D=Disgp\langle X\mid S\rangle:=[X^+]_\omega/\rho(S)
$$
for some generator set $X$ with defining relations $S\subseteq [X^+]_\omega\times[X^+]_\omega$, where $\rho(S)$ is the congruence on $([X^+]_\omega,\vdash,\dashv)$ generated by $S$.
For convenience, we let $[u]_m-[v]_n$ or $[u]_m=[v]_n$ stand for the pair $([u]_m,[v]_n)$ in $Disgp\langle X\rangle$.

\begin{lemm}(\cite[Theorem 5.1]{CZ})\label{tndsg}
Let $D=Disgp\langle X|S\rangle$ be the disemigroup generated
by $X$ with defining relations $S$ and $>$ be a monomial-center ordering on $[X^+]_\omega$.
If $S\subseteq Di\langle X\rangle$ is a Gr\"{o}bner-Shirshov basis,
then $Irr(S)=\{[u]_n\in [X^+]_\omega\mid [u]_n\neq \overline{[asb]_m} \ \mbox{ for any normal } S\mbox{-diword}\ [asb]_m \}$
is a set of normal forms of elements of $D$.
\end{lemm}

\section{Free digroup}

In this section, we first introduce the notion of a free digroup, and then provide a disemigroup presentation for it.
We will apply Lemma \ref{tndsg} to obtain the set of normal forms of elements of a free digroup which gives a construction of the free digroup.

Loday \cite{Lo99} used the term ``dimonoid" to refer to what we have called a disemigroup.
We have made a slight change in the terminology to be more consistent
with standard usage in semigroup theory. Let $(D,\vdash,\dashv)$ be a disemigroup. An element $e$ in  $D$ is called
a \emph{bar-unit} if it satisfies $e\vdash a = a\dashv e = a$ for all $a\in D$.
The following lemma shows that  Definition \ref{eq00} is equivalent to the Kinyon's Definition (\cite{Kinyon04}, Definition 4.1).

\begin{lemm}\label{l1}
Let $G$ be a disemigroup with a bar-unit $e$. Then $(G,e)$ is a digroup with the unit $e$ if
$$
\forall a\in G, \exists b\in G \ \mbox{such that} \ a\vdash b=b\dashv a=e.
$$
\end{lemm}
\noindent{\bf Proof.}
It suffices to prove that for all $a\in G$ there exists a unique  $b\in G$ such that $ a\vdash b=b\dashv a=e$.
Suppose that there exist $b$ and $b'$ such that
$
a\vdash b=b\dashv a=e, \ a\vdash b'=b'\dashv a=e.
$
Then $b=b\dashv e=b\dashv (a\vdash b')=(b\dashv a)\dashv b'=e\dashv b'$,
$b'=b'\dashv e=b'\dashv (a\vdash b)=(b'\dashv a)\dashv b=e\dashv b$.
It follows that
$$
b=e\dashv b'=e\dashv (e\dashv b)=(e\dashv e)\dashv b=e\dashv b=b'.  \ \ \square
$$

A map $f$ from a disemigroup $G$ to a disemigroup $H$ is called a \emph{disemigroup homomorphism}
if for any $a,b\in G$, $f(a\vdash b)=f(a)\vdash f(b),f(a\dashv b)=f(a)\dashv f(b)$.
Let $(G,1_G)$ and $(H,1_H)$ be digroups. A disemigroup homomorphism $f$ from $(G,1_G)$ to $(H,1_H)$ is called a \emph{digroup homomorphism}
if $f(1_G)=1_H$.

Note that if $f$ is  a digroup homomorphism from $(G,1_G)$ to $(H,1_H)$, then  for any $a\in G$, we have $f(a^{\dag})= f(a)^{\dag}$.

\begin{definition}\label{dfd}
Let $X$ be a subset of a digroup $(D,1_D)$. Then $(D,1_D)$ is called a \emph{free digroup} on $X$ provided the following holds:
for any digroup $(G,1_G)$ and any function $\varphi$ from the set $X$ into $(G,1_G)$,  there exists a unique extension of $\varphi$
to a digroup homomorphism $\varphi^*$ from $(D,1_D)$ into $(G,1_G)$.
\end{definition}

In the paper \cite{Ph}, it is shown that a free digroup on a set $X$ exists. Following from universal properties,
it is easy to see that the free digroup is uniquely determined up to isomorphism.
In this case, we denote the free digroup on $X$ by $(F(X),e)$.

We now turn to a construction of a free digroup.
Let a set $X$ be given. We write $X^{-1}$ which is a set disjoint from $X$ with a one-to-one correspondence, i.e.
$X^{-1}=\{x^{-1}\mid x\in X\}$ and $X\cap X^{-1}=\emptyset$,
and write $X^{\pm 1}=X\cup X^{-1}$.

The following lemma gives a disemigroup presentation to a free digroup on $X$.
Let $Y$ be a set, $S$ a set of relations in $Disgp\langle Y\rangle$ and $\rho(S)$ the congruence on~$Disgp\langle Y\rangle$ generated by~$S$.
For simplicity of notation, for all $[u]_m,[v]_n\in Disgp\langle Y\rangle$,
we write~$[u]_m=[v]_n$ in $Disgp\langle Y| S\rangle$ which means that~$([u]_m,[v]_n)\in \rho(S)$. We denote the equivalent
class of $[u]_m$ in $Disgp\langle Y\rangle$ by $[u]_m\rho(S)$.

\begin{lemm}\label{lemfdig}
Let $X$ be a set and $e\notin X^{\pm 1}$ be a symbol. Then $(Disgp\langle X^{\pm 1}\cup \{e\}| S\rangle,e\rho(S))$ is a free digroup on $X$,
where $S$ consists of the following relations:
$$
x\vdash x^{-1}=e=x^{-1}\dashv x,\  \ y\dashv e=y=e\vdash y,\ \ \ x\in X,\ y\in X^{\pm 1}\cup \{e\}.
$$
\end{lemm}

\noindent{\bf Proof.}
We first show that $(Disgp\langle X^{\pm 1}\cup \{e\}| S\rangle,e\rho(S))$ is a digroup. It is easy to see that $e\rho(S)$ is a bar-unit in $Disgp\langle X^{\pm 1}\cup \{e\}| S\rangle$.

Now,  in $Disgp\langle X^{\pm 1}\cup \{e\}| S\rangle$, the following holds:

(i) $y\vdash e=e\dashv y$ for all $y\in X^{\pm 1}\cup \{e\}$. In fact, $e\vdash e=e=e\dashv e$.
For any $x\in X$, $x\vdash e=x\vdash (x^{-1}\dashv x)=(x\vdash x^{-1})\dashv x=e\dashv x$.
For any $z^{-1}\in X^{-1}$, we have
$$
  z^{-1}\vdash e=z^{-1}\vdash (z\vdash z^{-1})=(z^{-1}\dashv z)\vdash z^{-1}\\=e\vdash z^{-1}=z^{-1}$$
  and
  $$
e\dashv z^{-1}=(z^{-1}\dashv z)\dashv z^{-1}=z^{-1}\dashv (z\vdash z^{-1})=z^{-1}\dashv e=z^{-1}.$$
It follows that $z^{-1}\vdash e=e\dashv z^{-1}$.

(ii) Define $x^0=e$ for all $x\in X$.
Then for any $[u]_m\in Disgp\langle X^{\pm 1}\cup \{e\}\rangle$,
we may assume that $u=x_1^{\varepsilon_1}\cdots x_n^{\varepsilon_n}$, $x_i\in X$, $\varepsilon_i\in \{-1,0,1\}$.
Then there exists $[x_n^{-\varepsilon_n}\cdots x_1^{-\varepsilon_1}e]_{n+1}$ such that
$$
[x_1^{\varepsilon_1}\cdots x_n^{\varepsilon_n}]_m \vdash [x_n^{-\varepsilon_n}\cdots x_1^{-\varepsilon_1}e]_{n+1}
=(x_1^{\varepsilon_1}\vdash \cdots\vdash x_n^{\varepsilon_n}\vdash x_n^{-\varepsilon_n} \vdash \cdots \vdash x_1^{-\varepsilon_1})\vdash e
=e\vdash e=e$$
and
\begin{eqnarray*}
\ [x_n^{-\varepsilon_n}\cdots x_1^{-\varepsilon_1}e]_{n+1}\dashv [x_1^{\varepsilon_1}\cdots x_n^{\varepsilon_n}]_m
&=&x_n^{-\varepsilon_n} \vdash \cdots \vdash x_1^{-\varepsilon_1}\vdash e \dashv x_1^{\varepsilon_1}\dashv \cdots\dashv x_n^{\varepsilon_n} \\
&=&e\dashv (x_n^{-\varepsilon_n} \dashv \cdots \dashv x_1^{-\varepsilon_1}\dashv x_1^{\varepsilon_1}\dashv \cdots\dashv x_n^{\varepsilon_n})
=e\dashv e=e.
\end{eqnarray*}
By Lemma \ref{l1}, $(Disgp\langle X^{\pm 1}\cup \{e\}| S\rangle,e\rho(S))$ is a digroup.

We now turn to show that $(Disgp\langle X^{\pm 1}\cup \{e\}| S\rangle,e\rho(S))$ is free. For any digroup $(G,1_G)$ and function $\varphi: X\rightarrow G$,
define a map
$$
\phi: X^{\pm 1}\cup \{e\}\rightarrow G, \ \ x\mapsto \varphi(x),\ x^{-1}\mapsto (\varphi(x))^\dag,\  e\rightarrow 1_G, \ x\in X.
$$
Then we obtain a disemigroup homomorphism
$$
\phi^*: Disgp\langle X^{\pm 1}\cup \{e\}\rangle \rightarrow G, \ \
[x_1^{\varepsilon_1}\cdots x_n^{\varepsilon_n}]_m\mapsto [\phi(x_1^{\varepsilon_1})\cdots \phi(x_n^{\varepsilon_n})]_m
$$
induced by $\phi$. Since
\begin{eqnarray*}
&&\phi^*(x\vdash x^{-1})=\varphi(x)\vdash (\varphi(x))^\dag=1_G=(\varphi(x))^\dag\dashv \varphi(x)=\phi^*(x^{-1}\dashv x), \ \ \phi^*(e)=1_G, \\
&&\phi^*(y\dashv e)=\phi^*(y)\dashv 1_G=\phi^*(y)=1_G \vdash \phi^*(y)=\phi^*(e\vdash y),
\end{eqnarray*}
we have a disemigroup homomorphism
$$
\varphi^*:Disgp\langle X^{\pm 1}\cup \{e\}| S\rangle\rightarrow G, \ \
[x_1^{\varepsilon_1}\cdots x_n^{\varepsilon_n}]_m\rho(S)\mapsto [\phi(x_1^{\varepsilon_1})\cdots \phi(x_n^{\varepsilon_n})]_m.
$$
As $\varphi^*(e\rho(S))=\phi(e)=1_G$, it follows that $\varphi^*$ is a digroup homomorphism.
It is easy to see that~$\varphi^*|_{X}=\varphi$  and such $\varphi^*$ is unique.
\ \ $\square$

\ \

In the following, we give an explicit construction of the free digroup~$(Disgp\langle X^{\pm 1}\cup \{e\}| S\rangle,e\rho(S))$
by using the Gr\"{o}bner-Shirshov bases method.

\begin{lemm}\label{lfreedi}
Let $X$ be a set and $e\notin X^{\pm 1}$ be a symbol.
Let
$$
S=\{[x x^{-1}]_2-e,[x^{-1} x]_1-e,[ye]_1-y,[ey]_2-y \mid x\in X, y\in X^{\pm 1}\cup \{e\}\}
$$
be a subset of $Di\langle X^{\pm 1}\cup \{e\}\rangle$. Then the following statements hold.
\begin{enumerate}
\item[(i)] \ $Disgp\langle X^{\pm 1}\cup \{e\} \mid S \rangle=Disgp\langle X^{\pm 1}\cup \{e\} \mid S\cup T \rangle$,
where $T$ consists of the following relations:
\begin{align*}
&t_1=[x^{-1} e]_2-x^{-1},\ \ \
t_2=[ex_{i_1}\cdots x_{i_n}x^{-1}]_1-[x_{i_1}\cdots x_{i_n} x^{-1}]_{n+1}, \\
&t_3=[x e]_2-[e x]_1, \ \ \ \ \
t_4=[x^{-1} x_{l_1}\cdots  x_{l_m} z^{-1}]_{m+2}-[x^{-1} x_{l_1}\cdots  x_{l_m} z^{-1}]_1,
\end{align*}
where $m\geq 0,\ n\geq 0, \ x,z,x_{l_k},x_{i_j}\in X,\ 1\leq k\leq m,\ 1\leq j\leq n$.
\item[(ii)] \ Let $X^{\pm 1}$ be a well-ordered set.
Define an ordering on $X^{\pm 1}\cup \{e\}$ by $e<z$ for
any $z\in X^{\pm 1}$ and the ordering on $[(X^{\pm 1}\cup \{e\})^+]_\omega$ is given in the
deg-lex-center way. Then $S\cup T$ is a Gr\"{o}bner-Shirshov basis in $Di\langle X^{\pm 1}\cup \{e\}\rangle$.
\end{enumerate}

\end{lemm}

\noindent{\bf Proof.}
$(i)$ Let $\rho(S)$ be the congruence of $Disgp\langle X^{\pm 1}\cup \{e\}\rangle$ generated by $S$. We only need to show that $T\subseteq \rho(S)$.
By the proof of Lemma \ref{lemfdig}, we obtain that the relations $t_1,t_3$ and $t_2$ with $n=0$ hold in $Disgp\langle X^{\pm 1}\cup \{e\} \mid S \rangle$.
In the disemigroup $Disgp\langle X^{\pm 1}\cup \{e\} \mid S \rangle$, we also have
\begin{align*}
(e\dashv x_{i_1}\dashv\cdots \dashv x_{i_n})\dashv x^{-1}&=(x_{i_1}\vdash\cdots \vdash x_{i_n}\vdash e)\dashv x^{-1}
=x_{i_1}\vdash\cdots \vdash x_{i_n}\vdash(e\dashv x^{-1}) \\
&=x_{i_1}\vdash\cdots \vdash x_{i_n}\vdash x^{-1} \ \ \ \mbox{if} \ n>0;\\
x^{-1}\vdash z^{-1}=x^{-1}\vdash (e\dashv z^{-1})&=(x^{-1}\vdash e)\dashv z^{-1}
=x^{-1}\dashv z^{-1};\\
x^{-1}\vdash x_{l_1}\vdash\cdots \vdash x_{l_m}\vdash z^{-1}&=x^{-1}\vdash x_{l_1}\vdash\cdots \vdash x_{l_m}\vdash (e\dashv z^{-1})
=x^{-1}\vdash (x_{l_1}\vdash\cdots \vdash x_{l_m}\vdash e)\dashv z^{-1} \\
&=x^{-1}\vdash (e\dashv x_{l_1}\dashv\cdots \dashv x_{l_m})\dashv z^{-1}
=(x^{-1}\vdash e)\dashv x_{l_1}\dashv\cdots \dashv x_{l_m}\dashv z^{-1} \\
&=x^{-1}\dashv x_{l_1}\dashv\cdots \dashv x_{l_m}\dashv z^{-1} \ \ \ \mbox{if} \ m>0.
\end{align*}
It follows that $t_2,t_4\in \rho(S)$.  Thus the result holds.

$(ii)$ Let $R=S\cup T$. We will show that all compositions in $R$ are trivial modulo $R$. Here we show that some compositions in $R$, as examples, are trivial modulo $R$.

It is evident that all possible compositions of left (right) multiplication are ones related to $t_4$ and they are equal to zero.
Let
$$
s_1=[x x^{-1}]_2-e,\ s_2=[x^{-1} x]_1-e,\ s_3=[y e]_1-y, \ s_4=[e y]_2-y, \ \ x\in X, y\in X^{\pm 1}\cup \{e\}.
$$

Here and subsequently, we denote by, for example, ``$f\wedge g,\ [w]_m$" the composition of the polynomials $f$ and  $g$ with ambiguity $[w]_m$.

1) Compositions of inclusion and left (right) multiplicative inclusion.

By noting that in $R$,
\begin{align*}
&s_3 \wedge s_4,\ w=ee,\ P(s_3)\cap P(s_4)=\emptyset; \ \ \ \ \ \
s_3 \wedge t_1,\ w=x^{-1}e,\ P(s_3)\cap P(t_1)=\emptyset;\\
&s_3 \wedge t_3,\ w=xe,\ P(s_3)\cap P(t_3)=\emptyset; \ \ \ \ \ \
s_4 \wedge t_2,\ w=ex^{-1},\ P(s_4)\cap P(t_2)=\emptyset; \\
&t_2 \wedge s_1,\ w=ex_{i_1}\cdots x_{i_n}x_{i_n}^{-1},\ P(t_5)\cap P([ex_{i_1}\cdots x_{i_{n-1}}s_1])=\{1\};\\
&t_2 \wedge s_4,\ w=ex_{i_1}\cdots x_{i_n}x^{-1},\ P(t_5)\cap P([s_4x_{i_2}\cdots x_{i_n}x^{-1}])=\emptyset; \\
&t_4 \wedge s_1,\ w=x^{-1}x_{l_1}\cdots x_{l_m}x_{l_m}^{-1},\ P(t_4)\cap P(x^{-1}x_{l_1}\cdots x_{l_{m-1}}s_1)=\{m+2\}; \\
&t_4 \wedge s_2,\ w=x_{l_1}^{-1}x_{l_1}\cdots x_{l_m}z^{-1},\ P(t_4)\cap P(s_2x_{l_2}\cdots x_{l_m}z^{-1})=\{m+2\},
\end{align*}
all possible of compositions of inclusion in $R$ are:
\begin{align*}
&t_2 \wedge s_1,\ [ex_{i_1}\cdots x_{i_n}x_{i_n}^{-1}]_1; \ \
t_4 \wedge s_1,\ [x^{-1}x_{l_1}\cdots x_{l_m}x_{l_m}^{-1}]_{m+2}; \ \
t_4 \wedge s_2,\ [x_{l_1}^{-1}x_{l_1}\cdots x_{l_m}z^{-1}]_{m+2},
\end{align*}
and all possible of compositions of left (right) multiplicative inclusion in $R$ are:
\begin{align*}
&s_3 \wedge s_4,\ [yee]_1,[eey]_3; \  \ \ \ \ \
s_3 \wedge t_1,\ [yx^{-1}e]_1,[x^{-1}ey]_3; \\
&s_3 \wedge t_3,\ [yxe]_1,[xey]_3; \ \ \ \ \ \
s_4 \wedge t_2,\ [yex^{-1}]_1,[ex^{-1}y]_3; \\
&t_2 \wedge s_4,\ [yex_{i_1}\cdots x_{i_n}x^{-1}]_1, [ex_{i_1}\cdots x_{i_n}x^{-1}y]_{n+3},
\end{align*}
where $m\geq 0, n\geq 0,\ y\in X^{\pm 1}\cup \{e\}, \ x,x_{l_k},x_{i_j},z\in X, 1\leq k\leq m, 1\leq j\leq n$.

It is easy to check that all above compositions are trivial modulo $R$. Here, for
example, we just check $t_4 \wedge s_1$ with ambiguity $[w]_{m+2}=[x^{-1}x_{l_1}\cdots x_{l_m}x_{l_m}^{-1}]_{m+2}$.
\begin{align*}
(t_4,s_1)_{[w]_{m+2}}&=x^{-1}\vdash x_{l_1}\vdash \cdots \vdash x_{l_{m-1}}\vdash e
-x^{-1}\dashv x_{l_1}\dashv \cdots \dashv x_{l_{m-1}}\dashv x_{l_m}\dashv x_{l_m}^{-1} \\
&\equiv x^{-1} \vdash e\dashv x_{l_1}\dashv \cdots \dashv x_{l_{m-1}}
-x^{-1}\dashv x_{l_1}\dashv \cdots \dashv x_{l_{m-1}}\dashv (x_{l_m}\vdash x_{l_m}^{-1}) \\
&\equiv x^{-1} \dashv x_{l_1}\dashv \cdots \dashv x_{l_{m-1}}-x^{-1}\dashv x_{l_1}\dashv \cdots \dashv x_{l_{m-1}}\dashv e
\equiv 0  \  \mathrm{mod} (R, [w]_{m+2}).
\end{align*}

2) Compositions of intersection and left (right) multiplicative intersection.

By noting that in $R$,
\begin{align*}
&s_1 \wedge s_2,\ w=xx^{-1}x,\ P([s_1x])\cap P([xs_2])=\{2\}; \\
&s_1 \wedge s_3,\ w=xx^{-1}e,\ P([s_1e])\cap P([xs_3])=\{2\};\\
&s_1 \wedge t_1,\ w=xx^{-1}e,\ P([s_1e])\cap P([xt_1])=\{3\}; \\
&s_1 \wedge t_4,\ w=xx^{-1}x_{l_1}\cdots x_{l_m}z^{-1},\ P([s_1x_{l_1}\cdots x_{l_m}z^{-1}])\cap P([xt_4])=\{m+3\}; \\
&s_2 \wedge s_1,\ w=x^{-1}xx^{-1},\ P([s_2x^{-1}])\cap P([x^{-1}s_1])=\{1,3\}; \\
&s_2 \wedge s_3,\ w=x^{-1}xe,\ P([s_2e])\cap P([x^{-1}s_3])=\{1\};\\
&s_2 \wedge t_3,\ w=x^{-1}xe,\ P([s_2e])\cap P([x^{-1}t_3])=\{1,3\}; \\
&s_3 \wedge s_3,\ w=yee,\ P([s_3e])\cap P([ys_3])=\{1\};\\
&s_3 \wedge s_4,\ w=yey',\ P([s_3y'])\cap P([ys_4])=\{1,3\}; \\
&s_3 \wedge t_2,\ w=yex_{i_1}\cdots x_{i_n}x^{-1},\ P([s_3x_{i_1}\cdots x_{i_n}x^{-1}])\cap P([yt_2])=\{1\};\\
&s_4 \wedge s_1,\ w=exx^{-1},\ P([s_4x^{-1}])\cap P([es_1])=\{3\}; \\
&s_4 \wedge s_2,\ w=ex^{-1}x,\ P([s_4x])\cap P([es_2])=\{2\}; \\
&s_4 \wedge s_3,\ w=eye,\ P([s_4e])\cap P([es_3])=\{2\}; \\
&s_4 \wedge s_4,\ w=eee,\ P([s_4e])\cap P([es_4])=\{3\}; \\
&s_4 \wedge t_1,\ w=ex^{-1}e,\ P([s_4e])\cap P([et_1])=\{3\}; \\
&s_4 \wedge t_2,\ w=eex_{i_1}\cdots x_{i_n}x^{-1},\ P([s_4x_{i_1}\cdots x_{i_n}x^{-1}])\cap P([et_2])=\{2\};\\
&s_4 \wedge t_3,\ w=exe,\ P([s_4e])\cap P([et_3])=\{3\}; \\
&s_4 \wedge t_4,\ w=ex^{-1}x_{l_1}\cdots x_{l_m}z^{-1},\ P([s_4x_{l_1}\cdots x_{l_m}z^{-1}])\cap P([et_4])=\{m+3\}; \\
&t_1 \wedge s_3,\ w=x^{-1}ee,\ P([t_1e])\cap P([x^{-1}s_3])=\{2\}; \\
&t_1 \wedge s_4,\ w=x^{-1}ey,\ P([t_1y])\cap P([x^{-1}s_4])=\{3\}; \\
&t_1 \wedge t_2,\ w=x^{-1}ex_{i_1}\cdots x_{i_n}z^{-1},\ P([t_1x_{i_1}\cdots x_{i_n}z^{-1}])\cap P([x^{-1}t_2])=\{2\};\\
&t_2 \wedge s_2,\ w=ex_{i_1}\cdots x_{i_n}x^{-1}x,\ P([t_2x])\cap P([ex_{i_1}\cdots x_{i_n}s_2])=\{1\}; \\
&t_2 \wedge s_3,\ w=ex_{i_1}\cdots x_{i_n}x^{-1}e,\ P([t_2e])\cap P([ex_{i_1}\cdots x_{i_n}s_3])=\{1\}; \\
&t_2 \wedge t_1,\ w=ex_{i_1}\cdots x_{i_n}x^{-1}e,\ P([t_2e])\cap P([ex_{i_1}\cdots x_{i_n}t_1])=\{1,n+3\}; \\
&t_2 \wedge t_4,\ w=ex_{i_1}\cdots x_{i_n}x^{-1}x_{l_1}\cdots x_{l_m}z^{-1},\
P([t_2x_{l_1}\cdots x_{l_m}z^{-1}])\cap P([ex_{i_1}\cdots x_{i_n}t_4])=\{n+m+3\};\\
&t_3 \wedge s_3,\ w=xee,\ P([t_3e])\cap P([xs_3])=\{2\}; \\
&t_3 \wedge s_4,\ w=xey,\ P([t_3y])\cap P([xs_4])=\{3\}; \\
&t_3 \wedge t_2,\ w=xex_{i_1}\cdots x_{i_n}z^{-1},\ P([t_3x_{i_1}\cdots x_{i_n}z^{-1}])\cap P([xt_2])=\{2\};\\
&t_4 \wedge s_2,\ w=x^{-1}x_{l_1}\cdots x_{l_m}z^{-1}z,\ P([t_4z])\cap P([x^{-1}x_{l_1}\cdots x_{l_m}s_2])=\{m+2\}; \\
&t_4 \wedge s_3,\ w=x^{-1}x_{l_1}\cdots x_{l_m}z^{-1}e,\ P([t_4e])\cap P([x^{-1}x_{l_1}\cdots x_{l_m}s_3])=\{m+2\}; \\
&t_4 \wedge t_1,\ w=x^{-1}x_{l_1}\cdots x_{l_m}z^{-1}e,\ P([t_4e])\cap P([x^{-1}x_{l_1}\cdots x_{l_m}t_1])=\emptyset; \\
&t_4 \wedge t_4,\ w=x^{-1}x_{l_1}\cdots x_{l_m}z^{-1}x_{i_1}\cdots x_{i_{n}}z'^{-1},\
P([t_4x_{i_1}\cdots x_{i_{n}}z'^{-1}])\cap P([x^{-1}x_{l_1}\cdots x_{l_m}t_4])=\emptyset,
\end{align*}
there is no composition of left (right) multiplicative intersection
and all possible of compositions of intersection in $R$ are:
\begin{align*}
s_1& \wedge s_2,\ [xx^{-1}x]_2; &
s_1& \wedge s_3,\ [xx^{-1}e]_2; &
s_1& \wedge t_1,\ [xx^{-1}e]_3; \\
s_1& \wedge t_4,\ [xx^{-1}x_{l_1}\cdots x_{l_m}z^{-1}]_{m+3}; &
s_2& \wedge s_1,\ [x^{-1}xx^{-1}]_1,[x^{-1}xx^{-1}]_3; &
s_2& \wedge s_3,\ [x^{-1}xe]_1;\\
s_2& \wedge t_3,\ [x^{-1}xe]_1,[x^{-1}xe]_3; &
s_3& \wedge s_3,\ [yee]_1;&
s_3& \wedge s_4,\ [yey']_1,[yey']_3; \\
s_3& \wedge t_2,\ [yex_{i_1}\cdots x_{i_n}x^{-1}]_1;&
s_4& \wedge s_1,\ [exx^{-1}]_3; &
s_4& \wedge s_2,\ [ex^{-1}x]_2; \\
s_4& \wedge s_3,\ [eye]_2; &
s_4& \wedge s_4,\ [eee]_3; &
s_4& \wedge t_1,\ [ex^{-1}e]_3; \\
s_4& \wedge t_2,\ [eex_{i_1}\cdots x_{i_n}x^{-1}]_2;&
s_4& \wedge t_3,\ [exe]_3; &
s_4& \wedge t_4,\ [ex^{-1}x_{l_1}\cdots x_{l_m}z^{-1}]_{m+3}; \\
t_1& \wedge s_3,\ [x^{-1}ee]_2; &
t_1& \wedge s_4,\ [x^{-1}ey]_3; &
t_1& \wedge t_2,\ [x^{-1}ex_{i_1}\cdots x_{i_n}z^{-1}]_2;\\
t_2& \wedge s_2,\ [ex_{i_1}\cdots x_{i_n}x^{-1}x]_1; &
t_2& \wedge s_3,\ [ex_{i_1}\cdots x_{i_n}x^{-1}e]_1; &
t_2& \wedge t_1,\ [ex_{i_1}\cdots x_{i_n}x^{-1}e]_1, \\
t_2& \wedge t_1,\ [ex_{i_1}\cdots x_{i_n}x^{-1}e]_{n+3}; &
t_2& \wedge t_4,\ [ex_{i_1}\cdots x_{i_n}x^{-1}x_{l_1}\cdots x_{l_m}z^{-1}]_{n+m+3};&
t_3& \wedge s_3,\ [xee]_2; \\
t_3& \wedge s_4,\ [xey]_3; &
t_3& \wedge t_2,\ [xex_{i_1}\cdots x_{i_n}z^{-1}]_2;&
t_4& \wedge s_2,\ [x^{-1}x_{l_1}\cdots x_{l_m}z^{-1}z]_{m+2}; \\
t_4& \wedge s_3,\ [x^{-1}x_{l_1}\cdots x_{l_m}z^{-1}e]_{m+2}, &
\end{align*}
where $m\geq 0,n\geq 0, \ y,y'\in X^{\pm 1}\cup \{e\},\ x,z,z',x_{l_k},x_{i_j}\in X, 1\leq k\leq m,1\leq j\leq n$.

It is easy to check that all compositions of intersection are trivial modulo $R$. Here, for
example, we just check $t_4 \wedge s_2$ with ambiguity $[w]_{m+2}=[x^{-1}x_{l_1}\cdots x_{l_m}z^{-1}z]_{m+2}$.
\begin{align*}
(t_4,s_2)_{[w]_{m+2}}&= x^{-1}\vdash x_{l_1}\vdash\cdots\vdash x_{l_m}\vdash e
-x^{-1}\dashv x_{l_1}\dashv\cdots\dashv x_{l_m}\dashv z^{-1}\dashv z \\
&\equiv x^{-1}\vdash e\dashv x_{l_1}\dashv\cdots\dashv x_{l_m}
-x^{-1}\dashv x_{l_1}\dashv\cdots\dashv x_{l_m}\dashv e \\
&\equiv x^{-1}\dashv x_{l_1}\dashv\cdots\dashv x_{l_m}
-x^{-1}\dashv x_{l_1}\dashv\cdots\dashv x_{l_m} \equiv 0 \ \mathrm{mod} (R, [w]_{m+2}).
\end{align*}

Therefore $R=S\cup T$ is a Gr\"{o}bner-Shirshov basis in $Di\langle X^{\pm 1}\cup \{e\}\rangle$.
\ \ $\square$

\ \

Let $X$ be a set and
$$
gp\langle X\rangle:=\{x_{i_1}^{\varepsilon_1}\cdots x_{i_n}^{\varepsilon_n}\mid x_{i_j}\in X, \varepsilon_j\in \{-1,1\},1\leq j\leq n, n\geq0
\ \mbox{ with } \varepsilon_j+\varepsilon_{j+1}\neq 0 \mbox{ whenever } x_{i_j}=x_{i_{j+1}}\}
$$
be the free group on $X$ with the unit $\epsilon$  (the empty word).
Then for any word $u\in (X^{\pm 1})^*$, there is a unique reduced word $\widehat{u}\in gp\langle X\rangle$.
Let $\tau:(X^{\pm 1}\cup \{e\})^*\rightarrow (X^{\pm 1})^*$ be the semigroup homomorphism which is deduced by
$$
x\mapsto x, \ \ x^{-1}\mapsto x^{-1}, \ \ e\mapsto \epsilon, \ \ \  x\in X.
$$
Let $v\in (X^{\pm 1}\cup \{e\})^*$. For abbreviation, we let $\widehat{v}$ stand for $\widehat{\tau(v)}$.

Let $u\in (X^{\pm1})^*$ and $u=x_{i_1}^{\varepsilon_1}x_{i_2}^{\varepsilon_2}\cdots x_{i_n}^{\varepsilon_n},n\geq 0$
with $x_{i_j}\in X$, $\varepsilon_j\in \{-1,1\}$.
Write $u^{-1}=x_{i_n}^{-\varepsilon_n}\cdots x_{i_2}^{-\varepsilon_2}x_{i_1}^{-\varepsilon_1}$ for $n>0$
and $\epsilon^{-1}=\epsilon$. Define
$$
\lambda: (X^{\pm 1})^*\backslash X^*\rightarrow \mathbb{N}, \ \
x_{i_1}^{\varepsilon_1}x_{i_2}^{\varepsilon_2}\cdots x_{i_n}^{\varepsilon_n}\mapsto \min\{j \mid \varepsilon_j=-1,1\leq j\leq n\}.
$$
For example, if $u=x_1x_2^{-1}x_3x_4^{-1}$, $x_i\in X$, then $\lambda(u)=2$.

The following theorem gives a characterization of the free digroup
$
(F(X),e)=(Disgp\langle X^{\pm 1}\cup \{e\}| S\rangle,e\rho(S)).
$
\begin{theorem}\label{freedigp}
Let $X$ be a set and $(F(X),e)$ the free digroup on $X$ with the unit $e$. Then
$$
F(X)=\Omega_e\cup\Omega_{X}\cup \Omega_{X^{-1}}
$$
where
\begin{align*}
\Omega_e&=\{[eu]_1\mid u\in X^*\},\\
\Omega_{X}&=\{[uxv]_{|u|+1} \mid x\in X, \ u,v\in gp\langle X\rangle\},\\
\Omega_{X^{-1}}&=\{[ux^{-1}v]_{|u|+1} \mid x\in X, \ u\in X^*,v\in (X^{\pm1})^*, \ ux^{-1}v\in gp\langle X\rangle\}.
\end{align*}
Moreover, in $F(X)$, the operations $\vdash$, $\dashv$ and $\dagger$ are as follows:
for any $[q]_n\in F(X), \ [p]_m\in \Omega_e\cup \Omega_{X^{-1}}, \ [uxv]_{|u|+1}\in \Omega_{X}$,
\begin{displaymath}
 [q]_n\vdash [p]_m=
\begin{cases}
\text{ $[e\widehat{qp}]_1$, if $\widehat{qp}\in X^*$,} \\
\text{ $[\widehat{qp}]_{_{\lambda(\widehat{qp})}}$, otherwise,}
\end{cases}
\ \ \ \ \ \ \ [p]_m\dashv [q]_n=
\begin{cases}
\text{ $[e\widehat{pq}]_1$, if $\widehat{pq}\in X^*$,} \\
\text{ $[\widehat{pq}]_{_{\lambda(\widehat{pq})}}$, otherwise,}
\end{cases}
\end{displaymath}
$$
[q]_n \vdash [uxv]_{|u|+1}=[\widehat{qu}xv]_{|\widehat{qu}|+1}, \ \ \ \ \ \ \ [uxv]_{|u|+1}\dashv [q]_n=[ux\widehat{vq}]_{|u|+1},
$$
\begin{displaymath}
  ([q]_n)^{\dagger}=
\begin{cases}
\text{ $[e(\widehat{q})^{^{-1}}]_1$, if $(\widehat{q})^{^{-1}}\in X^*$,} \\
\text{ $[(\widehat{q})^{^{-1}}]_{_{\lambda((\widehat{q})^{^{-1}})}}$, otherwise.}\ \ \ \ \ \ \ \ \ \ \ \ \ \ \ \ \ \ \ \ \ \ \ \ \ \ \ \ \ \ \ \ \ \ \ \ \ \ \ \ \
\end{cases}
\end{displaymath}

\end{theorem}
\noindent{\bf Proof.}
By Lemmas \ref{tndsg} and \ref{lfreedi}, we know that  $Irr(S\cup T)=\Omega_e\cup\Omega_{X}\cup \Omega_{X^{-1}}$ is a set of normal forms of elements of $F(X)=Disgp\langle X^{\pm 1}\cup \{e\} \mid S \rangle$. The operations formulas follow from relations in Lemma \ref{lfreedi}.
\ \ $\square$

\ \

Let $(D,1_D)$ be a digroup, $J$ be the set of all inverses in $D$ and $E$ be the halo of $D$, that is,
$$
J:=\{a^{\dag}\mid a\in D\}, \ E:=\{e\in D\mid e\vdash a=a\dashv e=a,  a\in D\}.
$$
Kinyon \cite{Kinyon04} showed that $J=\{a\vdash 1_D \mid a\in D\}$,
$J$ is a group in which $\vdash \ = \ \dashv$
and $(D,1_D)$ is isomorphic to the digroup $(E\times J,(1_D,1_D))$ where $\vdash$ and $\dashv$ are defined by
$$
(u,h)\vdash (v,k)=(h\vdash v\dashv h^{\dag},h \vdash k), \ \
(u,h)\dashv(v,k)=(u,h\dashv k).
$$
We shall call $J$ the \textit{group part} and $E$ the \textit{halo part} of the digroup $(D,1_D)$.

Now we consider  the group and halo parts of the free digroup $(F(X),e)$.
By Theorem \ref{freedigp},  we immediately have
\begin{coro}
Let $X$ be a set, $J$ and $E$ be the group and holo parts of $(F(X),e)$.
Then $J=\Omega_e\cup \Omega_{X^{-1}}$ and $E=\{e \}\cup \{[uxv]_{|u|+1}\in \Omega_X \mid \widehat{uxv}=\epsilon\ \mbox{ in }\ gp\langle X\rangle \}$.
\end{coro}

Let $X=\{x\}$. For simplicity of notation, we write $x^{-n}$  instead of $(x^{-1})^n$ for any $n\in \mathbb{N}$. Then for any $i\in \mathbb{Z}$, the length of $x^i$ is the absolute value $|i|$ of $i$. The following corollary gives the free digroup generated by a single element.

\begin{coro}
Let $X=\{x\}$.
Then  $F(X)=\Omega_e\cup\Omega_{x}\cup \Omega_{x^{-1}}$, where
$$
\Omega_e=\{[ex^n]_1\mid n\geq 0\},\ \
\Omega_{x}=\{[(x^ixx^j]_{|i|+1} \mid i,j\in \mathbb{Z} \},\ \
\Omega_{x^{-1}}=\{[x^{-m}]_1 \mid m\geq 1\},
$$
and the operations $\vdash$, $\dashv$ and $\dagger$ are as follows:
for any $[ex^n]_1,[ex^{n'}]_1\in \Omega_e, \ \ [x^ixx^j]_{|i|+1},[x^{i'}xx^{j'}]_{|i'|+1}\in \Omega_{x}$
and $[x^{-m}]_1,[x^{-m'}]_1\in \Omega_{x^{-1}}$, \ $t=i+1+j, \ t'=i'+1+j'$, \\
\begin{tabular}{ |c |c |c |c| }
\hline
\multirow{2}{1cm}{$\vdash$} &
\multirow{2}{2cm}{$[ex^{n'}]_1$}&
\multirow{2}{2cm}{$[x^{i'}xx^{j'}]_{|i'|+1}$}&
\multirow{2}{2cm}{$[x^{-m'}]_1$} \\
 & & & \\
\hline
\multirow{2}{2cm}{$[ex^{n}]_1$} &
\multirow{2}{2cm}{$[ex^{n+n'}]_1$} &
\multirow{2}{2cm}{$[x^{n+i'}xx^{j'}]_{|n+i'|+1}$} &
 $[ex^{n-m'}]_1 \ (n-m'\geq 0)$ \\
& &  &$[x^{n-m'}]_1 \ (n-m'<0)$ \\
\hline
\multirow{2}{2cm}{$[x^ixx^j]_{|i|+1}$} & $[ex^{t+n'}]_1 \ (t+n'\geq 0)$ &
\multirow{2}{2cm}{$[x^{t+i'}xx^{j'}]_{|t+i'|+1}$}
& $[ex^{t-m'}]_1 \ (t-m'\geq 0)$ \\
 & $[x^{t+n'}]_1 \ (t+n'< 0)$
 &  & $[x^{t-m'}]_1 \ (t-m'< 0)$ \\
\hline
\multirow{2}{2cm}{$[x^{-m}]_1$} &
$[ex^{-m+n'}]_1 \ (-m+n'\geq 0)$ &
\multirow{2}{4cm}{$[x^{-m+i'}xx^{j'}]_{|-m+i'|+1}$} &
\multirow{2}{2cm}{$[x^{-m-m'}]_1$}  \\
 & $[x^{-m+n'}]_1 \ (-m+n'<0)$ & &  \\
\hline
\end{tabular} \\
\begin{tabular}{ |c |c |c |c| }
\hline
\multirow{2}{1cm}{$\dashv$} &
\multirow{2}{2cm}{$[ex^{n'}]_1$}&
\multirow{2}{2cm}{$[x^{i'}xx^{j'}]_{|i'|+1}$}&
\multirow{2}{2cm}{$[x^{-m'}]_1$} \\
 & & & \\
\hline
\multirow{2}{2cm}{$[ex^{n}]_1$} &
\multirow{2}{2cm}{$[ex^{n+n'}]_1$} &
$[ex^{n+t'}]_1 \ (n+t'\geq 0)$ &
 $[ex^{n-m'}]_1 \ (n-m'\geq 0)$ \\
& &$[x^{n+t'}]_1 \ (n+t'< 0)$  &$[x^{n-m'}]_1 \ (n-m'<0)$ \\
\hline
\multirow{2}{2cm}{$[x^ixx^j]_{|i|+1}$} &
\multirow{2}{2cm}{$[x^ixx^{j+n'}]_{|i|+1}$} &
\multirow{2}{4cm}{$[x^{i}xx^{j+t'}]_{|i|+1}$} &
\multirow{2}{2cm}{$[x^{i}xx^{j-m'}]_{|i|+1}$} \\
 &  &  & \\
\hline
\multirow{2}{2cm}{$[x^{-m}]_1$} &
$[ex^{-m+n'}]_1 \ (-m+n'\geq 0)$ &
$[ex^{-m+t'}]_1 \ (-m+t'\geq 0)$ &
\multirow{2}{2cm}{$[x^{-m-m'}]_1$}  \\
 & $[x^{-m+n'}]_1 \ (-m+n'<0)$ &  $[x^{-m+t'}]_1 \ (-m+t'<0)$ &  \\
\hline
\end{tabular}

\begin{displaymath}
([ex^{n}]_1)^{\dag}=[x^{-n}]_1,  \ \ \  \ \
([x^{i}xx^{j}]_{|i|+1})^{\dag}=
\begin{cases}
\text{ $[ex^{-t}]_1$, if $t\leq 0$,} \\
\text{ $[x^{-t}]_1$, \ if $t>0$,}
\end{cases}
\ \ \ \ \  ([x^{-m}]_1)^{\dag}=[ex^{m}]_1.
\end{displaymath}

Moreover, $J=\{[ex^n]_1\mid n\geq 0\}\cup \{[x^{-m}]_1 \mid m\geq 1\}$ is the group part of $(F(X),e)$
and $E=\{e\}\cup \{[x^{-n}x^n]_{n+1} \mid n\geq 1\}\cup \{[x^nx^{-n}]_n \mid n\geq 1\}$ is the halo part of $(F(X),e)$.
\end{coro}

Let $(G,1_G)$ be a digroup and $H$ be a subset of $G$. A pair $(H,1_G)$ is called a \textit{subdigroup} of $(G,1_G)$
if $1_G\in H$ and $(H,1_G)$ is a digroup under the operations of $(G,1_G)$.
For a subset $U$ of $G$, we write $\langle U\rangle$ for the subdigroup of $G$ generated by $U$.

Let $(A,1_A)$ and $(B,1_B)$ be two digroups. The set of all digroup homomorphisms from $(A,1_A)$ to $(B,1_B)$ will be denoted by $Hom(A,B)$.
The set of all maps from a set $X$ to a set $Y$ will be denoted by $Y^X$ and the cardinality of the set $X$ will be denoted by $card(X)$.
In the following theorem, we are going to give a characterization of two isomorphic free digroups.

\begin{theorem}
Let $X$ and $Y$ be two sets. Then $(F(X),e)$ is isomorphic to $(F(Y),e)$ if and only if $card(X)=card(Y)$.
\end{theorem}
\noindent{\bf Proof.}
Following from universal properties of free digroups, it is easy to check that $(F(X),e)\cong (F(Y),e)$ if $card(X)=card(Y)$.
Let $\sigma$ be a digroup isomorphism from $(F(X),e)$ to $(F(Y),e)$ and $B$ be a digroup with $card(B)=2$.
Then $card(Hom(F(X),B))=card(Hom(F(Y),B))$. By Definition \ref{dfd}, we have $card(Hom(F(X),B))=card(B^X)=2^{card(X)}$
and $card(Hom(F(Y),B))=card(B^Y)=2^{card(Y)}$. It follows that $2^{card(X)}=2^{card(Y)}$.
If $X$ is finite, then so is $Y$ and $card(X)=card(Y)$. Otherwise, $Y$ is infinite and $card(X)\geq \aleph_0$.
For any $x\in X$, $\sigma(x)\in F(Y)$, there exists a finite subset $V_x$ of $Y$ such that
$\sigma(x)\in \langle V_x\rangle$. Let $W=\cup_{_{x\in X}}V_x$. Clearly $W\subseteq Y$. We claim that $W=Y$.
Otherwise, let $y\in Y\backslash W$. Then $\sigma^{-1}(y)\in F(X)$ and there exists a finite subset $U=\{x_1,\cdots, x_m\}$ of $X$
such that $\sigma^{-1}(y)\in \langle U\rangle$. Thus $y=\sigma(\sigma^{-1}(y))\in \sigma(\langle U\rangle)\subseteq \langle W\rangle$. But $y\in Y\backslash W$.
This contradicts the fact that $(F(Y),e)$ is the free digroup on $Y$. Therefore, $card(Y)=card(W)=card(\cup_{_{x\in X}}V_x)\leq card(X)\aleph_0=card(X)$.
The proof for $card(X)\leq card(Y)$ is similar and we have done.
\ \ $\square$

\ \

\noindent{\bf Acknowledgement}
We wish to express our thanks to the referee for helpful suggestions and comments.

\ \

\newcommand{\noopsort}[1]{}

\end{document}